\documentclass[a4paper,12pt]{article}
\usepackage{amsmath}
\usepackage{color}
\usepackage{amssymb}
\usepackage{latexsym}
\usepackage{url}
\numberwithin{equation}{section}

\def\R{{\bf R}}

\def\e{{\varepsilon}}

\newtheorem{thm}{Theorem}[section]

\newtheorem{lem}{Lemma}[section]

\newtheorem{rem}{Remark}[section]
\newtheorem{Def}{Definition}[section]

\title{Strauss exponent for semilinear wave equations with scattering space dependent damping}
\author{
Ning-An Lai
\footnote{Institute of Nonlinear Analysis and Department of Mathematics,
Lishui University,
Lishui 323000,
China.}
\footnote{
School of Mathematical Sciences, Fudan University, Shanghai 200433, China. e-mail: hyayue@gmail.com.
}
Ziheng Tu
\footnote{
 School of Data Science, Zhejiang University of Finance and Economics, 310018 Hangzhou, China. e-mail: tuziheng@zufe.edu.cn.
}
}

\date{
\[
\begin{array}{ll}
\mbox{\footnotesize{\bf Keywords:}}
& \mbox{\footnotesize semilinear wave equation, space dependent damping, blow-up, lifespan}\\
\mbox{\footnotesize{\bf MSC2010:}}
& \mbox{\footnotesize primary 35L71, secondary 35B44}\\
\end{array}
\]
}

\begin{document}
\maketitle
\begin{abstract}
It is believed or conjectured that the semilinear wave equations with scattering space dependent damping
admit the Strauss critical exponent, see Ikehata-Todorova-Yordanov \cite{ITY}(the bottom in page 2) and  Nishihara-Sobajima-Wakasugi \cite{N2}(conjecture iii in page 4). In this work, we are devoted to showing the conjecture is true at least when the decay rate of the space dependent variable coefficients before the damping is larger than 2. Also, if the nonlinear term depends only on the derivative of the solution, we may prove the upper bound of the lifespan is the same as that of the solution of the corresponding problem without damping. This shows in another way the \lq\lq hyperbolicity" of the equation.

\end{abstract}


\section{Introduction}
\par\quad
We consider the small data Cauchy problems
\begin{equation}
\label{Strauss}
\left\{
\begin{aligned}
& u_{tt} - \Delta u + \frac{\mu}{(1+|x|)^\beta} u_t = |u|^p, \quad\text{in $[0, T)\times\R^n$}, \\
& u(x,0)=\e f(x), \quad u_t(x,0)=\e g(x), \quad x\in\R^n,
\end{aligned}
\right.
\end{equation}
and
\begin{equation}
\label{Glassey}
\left\{
\begin{aligned}
& u_{tt} - \Delta u + \frac{\mu}{(1+|x|)^\beta} u_t = |u_t|^p, \quad\text{in $[0, T)\times\R^n$}, \\
& u(x,0)=\e f(x), \quad u_t(x,0)=\e g(x), \quad x\in\R^n,
\end{aligned}
\right.
\end{equation}
where $\mu>0, \beta>2$ are constants, the initial data $f(x), g(x)$ are compactly supported functions from the energy space
\[
f(x)\in H^1(\R^n),~~g(x)\in L^2(\R^n).
\]
Without loss of generality, we assume that
\begin{equation}\label{supp}
\begin{aligned}
supp~f(x), g(x)\in \{x:|x|\le 1\}.
\end{aligned}
\end{equation}

As mentioned in \cite{ITY}, equations \eqref{Strauss} and \eqref{Glassey} can be used to model the wave travel in a nonhomogeneous gas
with damping, and the space dependent coefficients represent the friction coefficients. Let us first take a look at the linear problem
\begin{equation}
\label{linear}
\left\{
\begin{aligned}
& u_{tt} - \Delta u + \frac{\mu}{(1+|x|)^\beta} u_t = 0, \quad\text{in $[0, T)\times\R^n$}, \\
& u(x,0)=\e f(x), \quad u_t(x,0)=\e g(x), \quad x\in\R^n.
\end{aligned}
\right.
\end{equation}
There are lots of literatures about the linear Cauchy problem \eqref{linear}. We list some but may be not all of them, i.e., \cite{IK1, IK2, ITY1, Mo, RTY1, RTY2, RTY3, RT, TY1, WY5}. Based on these known results, we may classify the linear problem \eqref{linear} into three cases, due to the value of decay rate $\beta$, see the table below.

\begin{center}
\begin{tabular}{|c|c|c|}
\hline
$\beta\in (-\infty, 1)$ & effective &
\begin{tabular}{c}
solution behaves like\\
that of heat equation
\end{tabular}
\\
\hline
$\beta=1$ &
\begin{tabular}{c}
scaling invariant\\
weak damping
\end{tabular} &
\begin{tabular}{c}
the asymptotic behavior\\depends on $\mu$
\end{tabular}
\\
 \hline
$\beta\in(1,\infty)$ & scattering &
\begin{tabular}{c}
solution behaves like
that\\ of wave equation without damping
\end{tabular}\\
\hline
\end{tabular}
\end{center}
We refer the reader to \cite{N2} for a good introduction to problem \eqref{linear}.

\begin{rem}
There is some little difference between the Cauchy problem \eqref{linear} and the correspond one with time dependent damping, i.e.,
\begin{equation}
\label{timedelinear}
\left\{
\begin{aligned}
& u_{tt} - \Delta u + \frac{\mu}{(1+t)^\beta} u_t = 0, \quad\text{in $[0, T)\times\R^n$}, \\
& u(x,0)=\e f(x), \quad u_t(x,0)=\e g(x), \quad x\in\R^n.
\end{aligned}
\right.
\end{equation}
Due to the results in \cite{HO, MN, N1, Wir1, Wir2, Wir3}, there is one more case for \eqref{timedelinear} than that of \eqref{linear} for different value of $\beta$, thus

\begin{center}
\begin{tabular}{|c|c|c|}
\hline
$\beta\in (-\infty, -1)$ & overdamping &
\begin{tabular}{c}
solution does not \\
decay to zero
\end{tabular}
\\
\hline
$\beta\in [-1, 1)$ & effective &
\begin{tabular}{c}
solution behaves like\\
that of heat equation
\end{tabular}
\\
\hline
$\beta=1$ &
\begin{tabular}{c}
scaling invariant\\
weak damping
\end{tabular} &
\begin{tabular}{c}
the asymptotic behavior\\depends on $\mu$
\end{tabular}
\\
 \hline
$\beta\in(1,\infty)$ & scattering &
\begin{tabular}{c}
solution behaves like
that\\ of wave equation without damping
\end{tabular}\\
\hline
\end{tabular}
\end{center}
We refer the reader to \cite{LT} for more detailed introduction.
\end{rem}

If we come back to the nonlinear problem \eqref{Strauss} and \eqref{Glassey}, as we consider small data problem, we always expect global existence for relatively big nonlinear power and blow-up in finite time for relatively small exponent. This means there is a boardline of value $p$ to distinguish these two cases, thus, the critical power. There are two kinds of critical powers related to the nonlinear problem \eqref{Strauss}, one is usually named with Strauss exponent($p_S(n)$), the critical power to the small data Cauchy problem of
\[
u_{tt}-\Delta u=|u|^p,
\]
while the other one is usually named with Fujita exponent($p_F(n)$), the critical power to the small data Cauchy problem of
\[
u_t-\Delta u=|u|^p.
\]
The Strauss exponent $p_S(n)$, which originates from the work \cite{Strauss}, is the positive root of the quadratic equation
\begin{equation}
\label{quadratic}
\gamma(p,n):=2+(n+1)p-(n-1)p^2=0,
\end{equation}
and $p_F(n)=1+\frac 2n$ is named after the pioneering work \cite{Fujita}. Noting that
\[
p_F(n)<p_S(n)=\frac{n+1+\sqrt{n^2+10n-7}}{2(n-1)}
\quad\mbox{for}\ n\ge2.
\]

For the semilinear problem \eqref{Strauss}, if the critical power equals to or at least is related to $p_S(n)$, we say the equation admits Strauss or \lq\lq wave" exponent, while if it has only connection to $p_F(n)$, we then say the problem admits Fujita or \lq \lq heat" exponent. Till the moment, most of the known results for the semilinear wave equations with space dependent damping focus on the effective and scaling invariant weak damping cases. Ikehata, Todorova and Yordanov \cite{ITY} proved that the critical power for \eqref{Strauss} with $0\le \beta<1$ is $p_C(n)=1+\frac{2}{n-\beta}$. Nishihara \cite{N3} also studied the same case but with absorbed semilinear term $|u|^{p-1}u$ and proved the diffusion phenomena. Li \cite{LiXin} considered the scaling invariant weak damping case($\beta=1$) with $\mu\ge n$ and proved the critical power is
\begin{equation}
\label{CriLi}
p_C(n)=\left\{
\begin{aligned}
& 1+\frac{2}{n-1},~~~n\ge 2, \\
& +\infty, ~~~~~~n=1.\\
\end{aligned}
\right.
\end{equation}
Obviously the above two cases are related to $p_F(n)$. Ikeda and Sobajima \cite{IS2} showed blow-up result for $\beta=1, n\ge 3, 0\le \mu<\frac{(n-1)^2}{n+1}$ and $\frac{n}{n-1}<p\le p_S(n+\mu)$. We still do not know what will happen for the gap $\mu\in [\frac{(n-1)^2}{n+1}, n)$ when $\beta=1$, this is the reason why we think the scaling invariant weak damping case is more delicate.
In \cite{N2} Nishihara, Sobajima and Wakasugi verified that the critical power is still $1+\frac{2}{n-\beta}$ when $\beta<0$. Actually, they studied the small data Cauchy problem \eqref{Strauss} with space-time dependent damping, but focusing on \lq\lq effective" case.

There are more results about the small data Cauchy problem of semilinear wave equations with time dependent damping, they are not completely solved, but at least for all the four cases of the corresponding linear problem \eqref{timedelinear} we have progresses.
For $\beta<-1$,  Ikeda and Wakasugi \cite{IWnew} verified the global existence result for all $p>1$. For $\beta\in [-1, 1)$, it has been proved that the critical is $p_F(n)$, see \cite{DLR13, FIW, ISW, LaZ, LZ, LNZ12, TY, WY17}. If $\beta=1$, we have two subcases. One is for relatively big $\mu$. D'Abbicco \cite{DABI}, D'Abbicco and Lucente \cite{DL1} and Wakasugi \cite{WY14_scale} have showed that
 the critical power is $p_F(n)$ when
 \[
 \mu\geq
 \left\{
 \begin{array}{cl}
5/3 & \mbox{for}\ n=1,\\
3 & \mbox{for}\ n=2,\\
n+2 & \mbox{for}\ n\geq 3.
\end{array}
\right.
\]
But if we fix $\mu=2$, then the results \cite{D-L, DLR15, KS, Lai, Pa, wak16} shows that the critical power is
\[
p_C(n)=\max\{p_F(n), p_S(n+2)\}.
\]
This reveals the fact that the critical power will move from $p_F(n)$ to $p_S(n)$ if $\mu$ is small enough, see \cite{LTW}. And the result of which was improved by \cite{IS} verifies the fact by showing blow-up result for
\[
0<\mu<\frac{n^2+n+2}{n+2}~and~1<p\le p_S(n+\mu).
\]
Unfortunately, still there is a gap for $\mu$ that we have no result. We come to the last case($\beta>1$). By introducing a bounded multiplier, the author and Takamura \cite{LT} proved the blow-up result for
\[
 1<p<
 \left\{
 \begin{array}{ll}
 p_S(n) & \mbox{for}\ n\ge2,\\
 \infty & \mbox{for}\ n=1.
 \end{array}
 \right.
 \]
Then Wakasa and Yordanov \cite{WY} proved non-global existence for $p=p_S(n)$. It seems that we can say $p=p_S(n)$ is exactly the critical power, due to the recent work by Liu and Wang \cite{LW}, in which they showed global existence for
\[
p>p_S(n)~and~n=3, 4.
\]

Problem \eqref{Glassey} is related to the so called Glassey conjecture, which was first proposed in \cite{Glassey}. It conjectures that the small data Cauchy problem of
\[
u_{tt}-\Delta u=a^2|\partial_tu|^p+b^2|\nabla u|^p,~~a^2+b^2>0
\]
admits a critical power
\[
p_G(n)=
 \left\{
 \begin{array}{ll}
 \frac{n+1}{n-1} & \mbox{for}\ n\ge2,\\
 \infty & \mbox{for}\ n=1.
 \end{array}
 \right.
\]
We expect that if we add scattering damping in the linear part, no matter it is time dependent or space dependent, the critical power will still be $p_G(n)$. It seems to be true at least for time dependent case, due to the work \cite{LT2}, in which the author and Takamura established blow-up and lifespan for
\[
1<p\le p_G(n),~~n\ge 1,
\]
and the work \cite{BL}, in which they proved global existence for $p>2, n=3$.

In this work, one of the goals is to show the blow-up power of the semilinear problem \eqref{Strauss} is related to $p_S(n)$,
when the decay rate of the space dependent coefficients before the damping is bigger than 2. The other goal is to establish blow-up result and lifespan estimate for the semilinear problem \eqref{Glassey} for $1<p\le p_G(n)$, with the same assumption on the space dependent coefficients. So far we know the efficient way to prove blow-up is the test function method, which was originated from the work of Zhang \cite{Z} for critical semilinear damped equation, and the work \cite{YZ} and \cite{Zhou} for semilinear wave equation. In \cite{Z} they used the smooth cut-off functions as the test functions, and it seems enough to get blow-up result for Fujita type power. But if we want to obtain blow-up result for Strauss type power, it is better to use some special solution of the linear wave equation as the test function, i.e.,
\[
F(t, x)=e^{-t}\int_{S^{n-1}}e^{x\cdot\omega}d\omega
 \]
 used in \cite{YZ} and the Gauss hypergeometric function used in \cite{Zhou} and \cite{Zhou2}. In \cite{IS} and \cite{IS2}, they succeeded in constructing Gauss hypergeometric function type solutions of the linear wave equations with time(space) dependent scaling invariant damping respectively, and then obtained blow-up results for Strauss type power for relatively small constant, by using a similar idea as that of \cite{Zhou2}. It seems that the scaling invariant damping has the same scale as that of the two terms of the linear wave equation, which makes it possible to construct the Gauss hypergeometric function type solutions. For the linear wave equation with scattering damping without such structure, we can not make it by the same way. However, we have a key observation: if we have $\phi(x)$ such that
\begin{equation}\label{test}
\Delta \phi-\frac{\mu}{(1+|x|)^{\beta}}\phi=\phi,
\end{equation}
which was first introduced in \cite{YZ06}, then $\Phi(t, x)=e^{-t}\phi(x)$ solves
\[
\partial_t^2\Phi-\Delta \Phi-\frac{\mu}{(1+|x|)^{\beta}}\partial_t\Phi=0.
\]
Based on this observation, we may borrow the idea of Ikeda-Sobajima-Wakasa \cite{ISWa} to establish desired blow-up results for \eqref{Strauss} and \eqref{Glassey} with $1<p\le p_S(n)$ and $1<p\le p_G(n)$ respectively.
\section{Main Result}
\par\quad
\begin{Def}\label{def1}
As in \cite{LT}, we say that $u$ is an energy solution of \eqref{Strauss} on $[0,T)$
if
\begin{equation}\nonumber
u\in C([0,T),H^1(\R^n))\cap C^1([0,T),L^2(\R^n))\cap L_{\rm loc}
^p(\R^n\times(0,T))
\end{equation}
satisfies $u(0, x)=\e f(x), u_t(0, x)=\e g(x)$ and
\begin{equation}\label{weaksol}
\begin{aligned}
&\e \int_{\R^n}g(x)\Psi(0, x)dx+\e\int_{\R^n}\frac{\mu}{(1+|x|)^{\beta}}f(x)\Psi(0, x)dx\\
&+\int_0^T\int_{\R^n}|u|^p\Psi(t, x)dxdt\\
=&-\int_0^T\int_{\R^n}u_t(t, x)\Psi_t(t, x)dxdt+\int_0^T\int_{\R^n}\nabla u(t, x)\cdot\nabla\Psi(t, x)dxdt\\
&-\int_0^T\int_{\R^n}\frac{\mu}{(1+|x|)^{\beta}}u(t, x)\Psi_t(t, x)dxdt\\
\end{aligned}
\end{equation}
for any $\Psi(t, x)\in C_0^{\infty}\left([0, T)\times \R^n\right)$.
\end{Def}
In the same way, we may define the energy solution to the Cauchy problem \eqref{Glassey}, the only thing we need to do is to replace the nonlinear term with $|u_t(t, x)|^p$.

\begin{Def}\label{def2}
We define the upper bound of lifespan to Problem \eqref{Strauss} or \eqref{Glassey} as
\[
T_{\e}=\sup \{T>0;~ there ~exists ~an ~energy~ solution ~to ~\eqref{Strauss}~ or~ \eqref{Glassey} ~in~ [0, T)\}.
\]
\end{Def}

The blow-up result for the Cauchy problem \eqref{Strauss} reads as
\begin{thm}
\label{thmStrauss}
Let $\beta>2$ and $1<p<p_S(n), n\ge 2$.
Assume that both $f\in H^1(\R^n)$ and $g\in L^2(\R^n)$ are non-negative and do not vanish identically. Also, the compact support assumption \eqref{supp} holds.
Suppose that an energy solution $u$ of \eqref{Strauss} satisfies
\begin{equation}
\label{support}
\mbox{\rm supp}\ u\ \subset\{(x,t)\in\R^n\times[0,T)\ :\ |x|\le t+1\}.
\end{equation}
Then there exists a constant $\e_0=\e_0(f,g,n,p,\mu, \beta)>0$
such that $T$ has to satisfy
\begin{equation}
\label{lifespan1}
T\le
 \left\{
 \begin{array}{ll}
 C\e^{-\frac{2(p-1)}{n+1-(n-1)p}} & \mbox{for}\ 1<p\le\frac{n}{n-1},\\
 C\varepsilon^{-2p(p-1)/\gamma(p, n)} & \mbox{for}\ \frac{n}{n-1}< p<p_S(n)
 \end{array}
 \right.
\end{equation}
for $0<\e\le\e_0$. Hereinafter, $C$ denotes a positive constant independent of $\e$ and may have different value from line to line.
\end{thm}

\begin{rem}
\label{remcritical1}
Compared to the classical upper bound of lifespan estimate for the semilinear wave equations without damping, the lifespan for $1<p\le \frac{n}{n-1}$ is not sharp, we leave it for further consideration.
\end{rem}

\begin{thm}
\label{thmStrausscri}
Let $\beta>2$ and $p=p_S(n), n\ge 2$.
Assume that both $f\in H^1(\R^n)$ and $g\in L^2(\R^n)$ are non-negative and do not vanish identically. Also, the compact support assumption \eqref{supp} holds.
Suppose that an energy solution $u$ of \eqref{Strauss} satisfies
\begin{equation}
\label{support}
\mbox{\rm supp}\ u\ \subset\{(x,t)\in\R^n\times[0,T)\ :\ |x|\le t+1\}.
\end{equation}
Then there exists a constant $\e_0=\e_0(f,g,n,p,\mu, \beta)>0$
such that $T$ has to satisfy
\begin{equation}
\label{lifespancri}
T\le \exp\left(C\e^{-p(p-1)}\right)
\end{equation}
for $0<\e\le\e_0$.
\end{thm}

\begin{thm}
\label{thmGlassey}
Let $\beta>2$ and $1<p \le p_G(n)$.
Assume that both $f\in H^1(\R^n)$ and $g\in L^2(\R^n)$ are non-negative and do not vanish identically. Also, the compact support assumption \eqref{supp} holds.
Suppose that an energy solution $u$ of \eqref{Glassey} satisfies
\begin{equation}
\label{support}
\mbox{\rm supp}\ u\ \subset\{(x,t)\in\R^n\times[0,T)\ :\ |x|\le t+1\}.
\end{equation}
Then there exists a constant $\e_0=\e_0(f,g,n,p,\mu, \beta)>0$
such that $T$ has to satisfy
\begin{equation}
\label{lifespan2}
T\le
 \left\{
 \begin{array}{ll}
 C\e^{-\left(\frac{1}{p-1}-\frac{n-1}{2}\right)^{-1}} & \mbox{for}\ 1<p<p_G(n),\\
 \exp\left(C\e^{-(p-1)}\right) & \mbox{for}\ p=p_G(n)
 \end{array}
 \right.
\end{equation}
for $0<\e\le\e_0$.
\end{thm}

\begin{rem}
In our main results, we have an assumption that $\beta> 2$. This restriction comes from the existence and asymptotic behavior of the solution $\phi(x)$ to equation \eqref{test}. Improvement from $\beta> 2$ to $\beta> 1$ is in our further consideration.
\end{rem}


\section{Proof for Theorem \ref{thmStrauss}}
\subsection{$\frac{n}{n-1}< p<p_S(n)$ part}
As we mentioned above, we will use the test function method similar as that of \cite{ISWa}, but we will use different test function.
Let $\eta(t)\in C^{\infty}([0, \infty))$ satisfy
\[
\eta(t)=
 \left\{
 \begin{array}{ll}
 1 & \mbox{for}\ t\le\frac12,\\
 decreasing & \mbox{for}\ \frac12<t<1,\\
 0 & \mbox{for}\ t\ge 1\\
 \end{array}
 \right.
\]
and
\[
|\eta'(t)|\le C,~|\eta''(t)|\le C.\\
\]
Let
\[
\eta_T(t)=\eta\left(\frac tT\right),~~~T\in (1, T_{\e}).
\]
We remark that one may assume $T_\e> 1$, otherwise our results hold obviously if we choose $\e$ to be small enough.
Choosing $\Psi=\eta_T^{2p'}(t)$ with $\frac1p+\frac{1}{p'}=1$ in \eqref{weaksol} and then integrating by parts one gets
\begin{equation}\label{thm1step1}
\begin{aligned}
&\e \int_{\R^n}g(x)dx+\e\int_{\R^n}\frac{\mu}{(1+|x|)^{\beta}}f(x)dx+\int_0^T\int_{\R^n}|u|^p\eta_T^{2p'}dxdt\\
=&-\int_0^T\int_{\R^n}\partial_t\left(u\partial_t\eta_T^{2p'}\right)dxdt+\int_0^T\int_{\R^n}u\partial_t^2\eta_T^{2p'}dxdt\\
&-\int_0^T\int_{\R^n}\frac{\mu}{(1+|x|)^{\beta}}u\partial_t\eta_T^{2p'}dxdt\\
=&\int_0^T\int_{\R^n}u\partial_t^2\eta_T^{2p'}dxdt-\int_0^T\int_{\R^n}\frac{\mu}{(1+|x|)^{\beta}}u\partial_t\eta_T^{2p'}dxdt\\
\triangleq&I_1+I_2.\\
\end{aligned}
\end{equation}
Noting that
\[
\begin{aligned}
&\partial_t\eta_T^{2p'}=\frac{2p'}{T}\eta_T^{2p'-1}\eta',\\
&\partial_t^2\eta_T^{2p'}=\frac{2p'(2p'-1)}{T^2}\eta_T^{2p'-2}|\eta'|^2+\frac{2p'}{T^2}\eta_T^{2p'-1}\eta'',\\
\end{aligned}
\]
we have
\begin{equation}\label{I}
\begin{aligned}
|I_1|&\le \frac{C}{T^2}\int_0^T\int_{\R^n}u\eta_T^{2p'-2}dxdt\\
&\le CT^{-2}\left(\int_0^T\int_{\R^n}|u|^p\eta_T^{2p'}dxdt\right)^{\frac1p}\left(\int_0^T\int_{\R^n}dxdt\right)^{\frac{1}{p'}}\\
&\le CT^{n+1-2p'}+\frac13\int_0^T\int_{\R^n}|u|^p\eta_T^{2p'}dxdt,\\
\end{aligned}
\end{equation}
where we have used the finite propagation speed property of the solution. In a similar way, we have
\begin{equation}\label{II}
\begin{aligned}
|I_2|
\le& CT^{-1}\left(\int_0^T\int_{\R^n}|u|^p\eta_T^{2p'}dxdt\right)^{\frac1p}\left(\int_0^T\int_{\R^n}\frac{1}{(1+|x|)^{p'
\beta}}dxdt\right)^{\frac{1}{p'}}\\
\le& CT^{-1}\left(\int_0^T\int_{\R^n}|u|^p\eta_T^{2p'}dxdt\right)^{\frac1p}\\
&\times\left(\int_0^T\int_{|x|\le t+1}\frac{(1+r)^{n-1-p'}}{(1+r)^{p'
(\beta-1)}}drdt\right)^{\frac{1}{p'}}\\
\le& CT^{n+1-2p'}+\frac13\int_0^T\int_{\R^n}|u|^p\eta_T^{2p'}dxdt.\\
\end{aligned}
\end{equation}
Here we should assume that
\[
n-p'>0\Rightarrow p> \frac{n}{n-1}.
\]
It holds by combining \eqref{thm1step1}, \eqref{I} and \eqref{II}
\begin{equation}\label{thm1step2}
\begin{aligned}
C_1(f, g)\e+\int_0^T\int_{\R^n}|u|^p\eta_T^{2p'}dxdt\le CT^{n-1-\frac{2}{p-1}},
\end{aligned}
\end{equation}
where
\[
C_1(f, g)=C\left(\int_{\R^n}g(x)dx+\int_{\R^n}\frac{\mu}{(1+|x|)^{\beta}}f(x)dx\right).
\]

\begin{lem}
If $\beta>0$, then for any $\alpha\in \R$ and a fixed constant $R$, we have
\begin{equation}\label{keyinq}
\begin{aligned}
\int_0^{t+R}(1+r)^{\alpha}e^{-\beta(t-r)}dr\le C(t+R)^{\alpha}.
\end{aligned}
\end{equation}
\end{lem}
{\bf Proof.} We split the proof into two parts. First it is easy to see
\begin{equation}\label{keyinq1}
\begin{aligned}
&\int_{\frac{t+R}{2}}^{t+R}(1+r)^{\alpha}e^{-\beta(t-r)}dr\\
\le& Ce^{-\beta t}(t+R)^{\alpha}\int_{\frac{t+R}{2}}^{t+R}e^{\beta r}dr\\
\le& C(t+R)^{\alpha}.\\
\end{aligned}
\end{equation}
On the other hand, we have
\begin{equation}\label{keyinq2}
\begin{aligned}
&\int_0^{\frac{t+R}{2}}(1+r)^{\alpha}e^{-\beta(t-r)}dr\\
\le& Ce^{-\frac{\beta t}{2}}\int_0^{\frac{t+R}{2}}(1+r)^{\alpha}dr\\
\le& Ce^{-\frac{\beta t}{2}}\times\left\{
 \begin{array}{ll}
 (t+R)^{\alpha+1} & \mbox{for}\ \alpha>-1\\
 \ln(t+R) & \mbox{for}\ \alpha=-1\\
 C & \mbox{for}\ \alpha<-1\\
 \end{array}
 \right.\\
 \le& C(t+R)^{\alpha}.\\
\end{aligned}
\end{equation}
Hence \eqref{keyinq} comes from \eqref{keyinq1} and \eqref{keyinq2}.

\begin{lem}\label{lemtest}
{\bf (Lemma 3.1 in \cite{YZ06})}. Assuming that $\beta>2$, then the following equation
\begin{equation}
\label{testequ}
\Delta \phi(x)-\frac{\mu}{(1+|x|)^\beta}\phi(x)=\phi(x),~~x\in \R^n
\end{equation}
admits a solution satisfying
\begin{equation}
\label{testpro}
\begin{aligned}
0<\phi(x)\le C(1+|x|)^{-\frac{n-1}{2}}e^{|x|}.
\end{aligned}
\end{equation}
\end{lem}

Let
\[
\Phi(t, x)=e^{-t}\phi(x),
\]
then it is easy to verify that
\begin{equation}
\label{testwave}
\begin{aligned}
\partial_t^2\Phi-\Delta \Phi-\frac{\mu}{(1+|x|)^{\beta}}\partial_t\Phi=0.
\end{aligned}
\end{equation}
Setting
\[
\Psi(t, x)= \eta_T^{2p'}(t)\Phi(t, x)
\]
in the definition of weak solution \eqref{weaksol} and making integration by parts, we get
\begin{equation}\label{thm1step3}
\begin{aligned}
&\e \int_{\R^n}g(x)\phi(x)dx+\e\int_{\R^n}\left(1+\frac{\mu}{(1+|x|)^{\beta}}\right)f(x)\phi(x)dx\\
&+\int_0^T\int_{\R^n}|u|^p\eta_T^{2p'}\Phi dxdt\\
=&\int_0^T\int_{\R^n}u\left(\partial_t^2\eta_T^{2p'}\Phi+2\partial_t\eta_T^{2p'}\partial_t\Phi-\frac{\mu}{(1+|x|)^{\beta}}
\partial_t\eta_T^{2p'}\Phi\right)dxdt\\
&+\int_0^T\int_{\R^n}u\eta_T^{2p'}\left(\partial_t^2\Phi-\Delta\Phi-\frac{\mu}{(1+|x|)^{\beta}}\partial_t\Phi\right)dxdt\\
=&\int_0^T\int_{\R^n}u\left(\partial_t^2\eta_T^{2p'}\Phi+2\partial_t\eta_T^{2p'}\partial_t\Phi-\frac{\mu}{(1+|x|)^{\beta}}
\partial_t\eta_T^{2p'}\Phi\right)dxdt\\
\triangleq&I_3+I_4+I_5.\\
\end{aligned}
\end{equation}
Noting that all the three remaining terms in the last equality include at least one derivative of the smooth cut-off function $\eta_T(t)$, this fact is the key to obtain our desired result. By the similar way as \eqref{I}, combining \eqref{testpro} and \eqref{keyinq},
we can estimate $I_3, I_4, I_5$ as
\begin{equation}\label{III}
\begin{aligned}
|I_3|\le& \int_0^T\int_{\R^n}|u\partial_t^2\eta_T^{2p'}\Phi|dxdt\\
\le& CT^{-2}\left(\int_0^T\int_{\R^n}\eta_T^{2p'}|u|^p dxdt\right)^{\frac1p}\left(\int_0^T\int_{\R^n}\Phi^{p'}dxdt\right)^{\frac{1}{p'}}\\
\le& CT^{-2}\left(\int_0^T\int_{\R^n}\eta_T^{2p'}|u|^p dxdt\right)^{\frac1p}\\
&\times\left(\int_{\frac T2}^T\int_{r\le t+1}e^{-p'(t-r)}(1+r)^{n-1-\frac{n-1}{2}p'}dxdt\right)^{\frac{1}{p'}}\\
\le& CT^{-2+(n-\frac{n-1}{2}p')\frac{1}{p'}}\left(\int_0^T\int_{\R^n}\eta_T^{2p'}|u|^p dxdt\right)^{\frac1p}.\\
\end{aligned}
\end{equation}
\begin{equation}\label{V}
\begin{aligned}
|I_4|\le& \int_0^T\int_{\R^n}|u\partial_t\eta_T^{2p'}\partial_t\Phi|dxdt\\
\le& CT^{-1}\left(\int_0^T\int_{\R^n}\eta_T^{2p'}|u|^p dxdt\right)^{\frac1p}\left(\int_0^T\int_{\R^n}\Phi^{p'}dxdt\right)^{\frac{1}{p'}}\\
\le& CT^{-1+(n-\frac{n-1}{2}p')\frac{1}{p'}}\left(\int_0^T\int_{\R^n}\eta_T^{2p'}|u|^p dxdt\right)^{\frac1p}\\
\end{aligned}
\end{equation}
and
\begin{equation}\label{IV}
\begin{aligned}
|I_5|\le C|I_4|\le CT^{-1+(n-\frac{n-1}{2}p')\frac{1}{p'}}\left(\int_0^T\int_{\R^n}\eta_T^{2p'}|u|^p dxdt\right)^{\frac1p}.\\
\end{aligned}
\end{equation}
We then conclude from \eqref{thm1step3}, \eqref{III}, \eqref{V} and \eqref{IV} that
\begin{equation}\label{thm1step4}
\begin{aligned}
C_2(f, g)\e\le CT^{-1+(n-\frac{n-1}{2}p')\frac{1}{p'}}\left(\int_0^T\int_{\R^n}\eta_T^{2p'}|u|^p dxdt\right)^{\frac1p},\\
\end{aligned}
\end{equation}
where
\[
C_2(f, g)=\int_{\R^n}g(x)\phi(x)dx+\int_{\R^n}\left(1+\frac{\mu}{(1+|x|)^{\beta}}\right)f(x)\phi(x)dx.\\
\]
This in turn implies that
\begin{equation}\label{thm1step5}
\begin{aligned}
\left(C_2(f, g)\e\right)^pT^{n-\frac{n-1}{2}p}\le\int_0^T\int_{\R^n}\eta_T^{2p'}|u|^p dxdt.\\
\end{aligned}
\end{equation}
Noting the assumption \eqref{supp}, $C_2(f, g)$ is nothing but a finite constant, then \eqref{thm1step2} and \eqref{thm1step5} yield
\[
T\le C\e^{-\frac{2p(p-1)}{\gamma(p, n)}},
\]
which is the second part of lifespan estimate \eqref{lifespan1}.
\subsection{$1<p\le\frac{n}{n-1}$ part}

In this subsection, we aim to study the upper bound of lifespan estimate for relatively small powers, i.e., $1<p\le\frac{n}{n-1}$. We have showed in \eqref{thm1step3}
\begin{equation}\label{thm1step6}
\begin{aligned}
&\e \int_{\R^n}g(x)\phi(x)dx+\e\int_{\R^n}\left(1+\frac{\mu}{(1+|x|)^{\beta}}\right)f(x)\phi(x)dx\\
&+\int_0^T\int_{\R^n}|u|^p\eta_T^{2p'}\Phi dxdt\\
=&\int_0^T\int_{\R^n}u\left(\partial_t^2\eta_T^{2p'}\Phi+2\partial_t\eta_T^{2p'}\partial_t\Phi-\frac{\mu}{(1+|x|)^{\beta}}
\partial_t\eta_T^{2p'}\Phi\right)dxdt\\
\triangleq&III+V+IV.\\
\end{aligned}
\end{equation}
We will re-estimate the three terms on the right hand side as follows.
\begin{equation}\label{III1}
\begin{aligned}
|I_3|\le& \int_0^T\int_{\R^n}|u\partial_t^2\eta_T^{2p'}\Phi|dxdt\\
\le& CT^{-2}\left(\int_0^T\int_{\R^n}\eta_T^{2p'}\Phi|u|^p dxdt\right)^{\frac1p}\left(\int_0^T\int_{\R^n}\Phi dxdt\right)^{\frac{1}{p'}}\\
\le& CT^{-2}\left(\int_0^T\int_{\R^n}\eta_T^{2p'}\Phi|u|^p dxdt\right)^{\frac1p}\\
&\times\left(\int_{\frac T2}^T\int_{r\le t+1}e^{-(t-r)}(1+r)^{n-1-\frac{n-1}{2}}dxdt\right)^{\frac{1}{p'}}\\
\le& CT^{-2+\frac{n+1}{2p'}}\left(\int_0^T\int_{\R^n}\eta_T^{2p'}\Phi|u|^p dxdt\right)^{\frac1p}.\\
\end{aligned}
\end{equation}
\begin{equation}\label{V1}
\begin{aligned}
|I_4|\le& \int_0^T\int_{\R^n}|u\partial_t\eta_T^{2p'}\partial_t\Phi|dxdt\\
\le& CT^{-1}\left(\int_0^T\int_{\R^n}\eta_T^{2p'}\Phi|u|^p dxdt\right)^{\frac1p}\left(\int_0^T\int_{\R^n}\Phi dxdt\right)^{\frac{1}{p'}}\\
\le& CT^{-1+\frac{n+1}{2p'}}\left(\int_0^T\int_{\R^n}\eta_T^{2p'}\Phi|u|^p dxdt\right)^{\frac1p}\\
\end{aligned}
\end{equation}
and
\begin{equation}\label{IV1}
\begin{aligned}
|I_5|\le C|I_4|\le CT^{-1+\frac{n+1}{2p'}}\left(\int_0^T\int_{\R^n}\eta_T^{2p'}\Phi|u|^p dxdt\right)^{\frac1p}.\\
\end{aligned}
\end{equation}
Hence we have by combining \eqref{thm1step6}, \eqref{III1}, \eqref{V1} and \eqref{IV1} that
\begin{equation}\label{thm1step7}
\begin{aligned}
&\e \int_{\R^n}g(x)\phi(x)dx+\e\int_{\R^n}\left(1+\frac{\mu}{(1+|x|)^{\beta}}\right)f(x)\phi(x)dx\\
&+\int_0^T\int_{\R^n}|u|^p\eta_T^{2p'}\Phi dxdt\\
\le &CT^{-1+\frac{n+1}{2p'}}\left(\int_0^T\int_{\R^n}\eta_T^{2p'}\Phi|u|^p dxdt\right)^{\frac1p}\\
\le &CT^{-p'+\frac{n+1}{2}}+\frac12\int_0^T\int_{\R^n}\eta_T^{2p'}\Phi|u|^p dxdt,\\
\end{aligned}
\end{equation}
this yields
\begin{equation}\label{thm1step8}
\begin{aligned}
T\le C\e^{-\frac{2(p-1)}{n+1-(n-1)p}},~~for~1<p<\frac{n+1}{n-1},\\
\end{aligned}
\end{equation}
which is the first part of the lifespan estimate \eqref{lifespan1} for $1<p\le\frac{n}{n-1}$. And we finish the proof of Theorem \ref{thmStrauss}.

\begin{rem}
For the power
\[
p\in \left(\frac{n}{n-1}, \frac{n+1}{n-1}\right),
\]
we will use the second part of the lifespan estimate \eqref{lifespan1}, due to the reason
\[
\frac{2p(p-1)}{\gamma(p, n)}<\frac{2(p-1)}{n+1-(n-1)p}.\\
\]
\end{rem}

\section{Proof for Theorem \ref{thmStrausscri}}

In this section, we focus on the \lq\lq critical" case $p=p_s(n)$. Again we will use the idea in \cite{ISWa}.
However, as the wave-scaling property is loss in our model, we could not expect to find the hypergeometric function as the exact test function. We shall apply the argument in \cite{TL} to construct the test function, i.e.,
\begin{equation}\label{bq}
b_q(t, x)=\int_0^1e^{-\eta t}\psi_{\eta}(x)\eta^{q-1}d\eta,\ \ \ q>0
\end{equation}
where $\psi_\eta(x)$ is the solution of
$$\left(\Delta-\frac{\mu\eta}{(1+|x|)^\beta}\right)\psi_\eta(x)=\eta^2\psi_\eta(x).$$

\begin{lem}\label{testcri}
Assuming
\[
V(x)=\frac{\mu}{(1+|x|)^\beta}
\]
and $\beta>2$. Then for given $\eta\in [0,1]$, there exists function
 $\psi_\eta\in C^2(\mathbb{R}^n)$ satisfying
\begin{equation}\label{testequcri}
\Delta \psi_\eta-\eta V\psi_\eta=\eta^2\psi_\eta
\end{equation}
such that for large $|\eta x|$,
\begin{equation}\label{asy_psi}
\psi_\eta(x)\sim \varphi_\eta(x):=\int_{\mathbb{S}^{n-1}}e^{\eta x\omega}d\omega(\sim|\eta x|^{\frac{1-n}2}e^{|\eta x|}).
\end{equation}
\end{lem}
The proof of the above lemma is parallel to that of Lemma 3.1 in \cite{YZ06}, and we postpone it to the appendix.

With $\psi_\eta(x)$ in hand, We have for $b_q(x,t)$
\begin{lem}\label{lembq}
$(i)$  $b_{q}(t, x)$ satisfies following identities
$$\frac{\partial}{\partial t}b_{q}(t, x)=-b_{q+1}(t, x),\ \ \frac{\partial^2}{\partial t^2}b_{q}(t, x)=b_{q+2}(t, x)$$and
\begin{equation}\label{bq2}
\Delta b_q(t, x)=V\cdot b_{q+1}(t, x)+b_{q+2}(t, x),
\end{equation}
and hence, $b_q(t, x)$ satisfies
\begin{equation}\label{1bq2}
\partial^2_tb_q-\Delta b_q-V\partial_tb_q=0.
\end{equation}
$(ii)$ For any positive constant $R>1$, $b_q(t, x)$ satisfies the following asymptotic behavior
\begin{equation}\label{bq3}
b_q(t, x)\sim\left\{
\begin{array}{ll}
(t+R+|x|)^{-q}\ &\mbox{if}\ 0<q<\frac{n-1}2,\\
(t+R+|x|)^{\frac{-n+1}2}(t+R-|x|)^{\frac{n-1}{2}-q}\ &\mbox{if}\ q>\frac{n-1}2.\\
\end{array}
\right.
\end{equation}
\end{lem}
{\bf proof.}
$(i)$ These identities can be proved by direct computation, We omit the details.

$(ii)$
Due to \eqref{asy_psi}, we have for any positive constant $R>1$,
$$b_q(t, x)=\int_0^1e^{-\eta t}\psi_{\eta}(x)\eta^{q-1}d\eta\sim\int_0^1e^{-\eta (t+R)}\varphi_{\eta}(x)\eta^{q-1}d\eta.$$
Further applying the plane wave formula to $\varphi_{\eta}(x)$ (see \cite{FJohn}, page 8), we have
\begin{eqnarray*}
b_q(t, x)&\sim&\int_0^1e^{-\eta (t+R)}\omega_{n-1}\int_{-1}^1(1-\theta^2)^{\frac{n-3}{2}}e^{\theta\eta|x|}d\theta\eta^{q-1}d\eta\\
&=&\omega_{n-1}\int_{-1}^1(1-\theta^2)^{\frac{n-3}2}d\theta\int_0^1e^{\eta(\theta|x|-t-R)}\eta^{q-1}d\eta\\
&=&\omega_{n-1}\int_{-1}^1(1-\theta^2)^{\frac{n-3}2}(t+R-\theta|x|)^{-q}d\theta\int_0^{t+R-\theta|x|}e^{-\zeta}\zeta^{q-1}d\zeta
\end{eqnarray*}
where $\zeta=(t+R-\theta|x|)\eta$ and $q>0$. Noting that for $-1\leqslant\theta\leqslant1$, we have $R-1\leqslant t+R-\theta|x|<\infty$, and hence
$$c_0\leqslant\int_0^{t+R-\theta|x|} e^{-\zeta}\zeta^{q-1}d\zeta\leqslant C_0,$$
with $c_0=\int_0^{R-1} e^{-\zeta}\zeta^{q-1}d\zeta$ and $C_0=\int_0^\infty e^{-\zeta}\zeta^{q-1}d\zeta=\Gamma(q)$.
Thus
\begin{equation}\label{hyper}
\begin{aligned}
b_q(t, x)&\sim\int_{-1}^1(1-\theta^2)^{\frac{n-3}2}(t+R-\theta|x|)^{-q}d\theta\\
&=2^{n-2}\int_{0}^1(\tilde{\theta}(1-\tilde{\theta}))^{\frac{n-3}2}(t+R+|x|-2\tilde{\theta}|x|)^{-q}d\tilde{\theta}\\
&=2^{n-2}(t+R+|x|)^{-q}\int_{0}^1(\tilde{\theta}(1-\tilde{\theta}))^{\frac{n-3}2}(1-\tilde{\theta}z)^{-q}d\tilde{\theta}.\\
\end{aligned}
\end{equation}
Here $\tilde{\theta}=\frac{\theta+1}2$ and $z=\frac{2|x|}{t+R+|x|}$. Direct analysis on integral shows for $n\geqslant 2$,
$$h(z)=\int_{0}^1(\tilde{\theta}(1-\tilde{\theta}))^{\frac{n-3}2}(1-\tilde{\theta}z)^{-q}d\tilde{\theta}$$
is integrable for $0\leqslant z<1$. Since for $z=1$, if $0<q<\frac{n-1}2$, $h(z)$ is continuous at $1$, thus bounded over $[0,1]$. Otherwise $q>\frac{n-1}2$,
$h(z)$ behaves as $(1-z)^{\frac{n-1}2-q}$ around $z=1$. It then follows the asymptotic behavior \eqref{bq3}.

\begin{rem}
In fact, for $\gamma>\beta>0$, the hypergeometric function has following integral representation
$$F(\alpha,\beta,\gamma;z)=\frac{\Gamma(\gamma)}{\Gamma(\beta)\Gamma(\gamma-\beta)}\int_0^1t^{\beta-1}(1-t)^{\gamma-\beta-1}(1-zt)^{-\alpha}dt,\ |z|<1.$$
Then \eqref{hyper} gives that
\begin{equation*}
b_q(t, x)\sim(t+R+|x|)^{-q}F(q,\frac{n-1}2,n-1;\frac{2|x|}{t+R+|x|}),
\end{equation*}
which means that our test function actually has the same asymptotic property as the hypergeometric function used in \cite{Zhou2}.
\end{rem}

As in \cite{ISWa}, we introduce
\[
\theta(t)=
 \left\{
 \begin{array}{ll}
 0 & \mbox{for}\ t<\frac12,\\
 \eta(t) & \mbox{for}\ t\ge\frac12,\\
 \end{array}
 \right.
{\ \ \ \ \ \theta_{M}(t):=\theta(\frac{t}{M}).}
\]
For $M\in (1, T)$, utilizing Lemma 3.1 in \cite{ISWa}, it concludes from \eqref{thm1step5} and the steps to get it
\begin{equation}\label{thmcristep1}
\begin{aligned}
\int_0^T\int_{\R^n}\theta_M^{2p'}|u|^pdxdt\ge\int_0^M\int_{\R^n}\theta_M^{2p'}|u|^pdxdt\ge C\e^pM^{n-\frac{n-1}{2}p}.
\end{aligned}
\end{equation}
Set
\[
Y[w](M)=\int_1^M\left(\int_0^T\int_{\R^n}w(t, x)\theta_\sigma^{2p'}(t)dxdt\right)\sigma^{-1}d\sigma,
\]
then for $q=\frac{n-1}{2}-\frac1p$ we have by combining \eqref{bq3} and \eqref{thmcristep1}
\begin{equation}\label{Ycrideriv}
\begin{aligned}
M\frac{d}{dM}Y\left[|u|^pb_q(t, x)\right](M)=&\int_0^T\int_{\R^n}\theta_M^{2p'}b_q|u|^pdxdt\\
\ge& CM^{-\left(\frac{n-1}{2}-\frac1p\right)}\int_0^T\int_{\R^n}\theta_M^{2p'}|u|^pdxdt\\
\ge& C\e^p,
\end{aligned}
\end{equation}
where we used the fact
\[
n-\frac{n-1}{2}p=\frac{n-1}{2}-\frac1p
\]
for $p=p_S(n)$. Also, direct computation shows that
\begin{equation}\label{Ycri}
\begin{aligned}
Y\left[b_q|u|^p\right](M)
=&\int_1^M
\left(\int_0^T\int_{\R^n}b_q|u|^p\theta_\sigma^{2p'}(t)dxdt\right)\sigma^{-1}d\sigma\\
=&\int_0^T\int_{\R^n}b_q|u|^p\int_1^M\theta_\sigma^{2p'}\sigma^{-1}d\sigma dxdt\\
=&\int_0^T\int_{\R^n}b_q|u|^p\int_{\frac tM}^t\theta^{2p'}(s)s^{-1}ds dxdt\\
\le&\int_0^{\frac M2}\int_{\R^n}b_q|u|^p\eta^{2p'}\left(\frac tM\right)\int_{\frac 12}^1 s^{-1}ds dxdt\\
&+\int_{\frac M2}^T\int_{\R^n}b_q|u|^p\int_{\frac tM}^1 \theta^{2p'}(s)s^{-1}ds dxdt\\
\le&\int_0^{\frac M2}\int_{\R^n}b_q|u|^p\eta^{2p'}\left(\frac tM\right)\int_{\frac 12}^1 s^{-1}ds dxdt\\
&+\int_{\frac M2}^T\int_{\R^n}b_q|u|^p\theta^{2p'}\left(\frac tM\right)\int_{\frac tM}^1 s^{-1}ds dxdt\\
\le&\int_{0}^T\int_{\R^n}b_q|u|^p\eta^{2p'}\left(\frac tM\right)\int_{\frac 12}^1 s^{-1}ds dxdt\\
\le& C\log2\int_0^T\int_{\R^n}\eta_M^{2p'}b_q|u|^pdxdt.\\
\end{aligned}
\end{equation}
For simplicity, we denote $Y(M)$ for $Y\left[b_q|u|^p\right](M)$, then \eqref{Ycri} yields
\begin{equation}\label{Ycriineq}
\begin{aligned}
Y(M)\le& C\int_0^T\int_{\R^n}\eta_M^{2p'}b_q|u|^pdxdt\\
=&C\int_0^T\int_{\R^n}\left(\partial_t^2u-\Delta u+\frac{\mu}{(1+|x|)^\beta}\partial_tu\right)b_q\eta_{M}^{2p'}\\
=&-C\Big(\e\int_{\R^n}g(x)b_q(0, x)dx+\e\int_{\R^n}f(x)b_{q+1}(0, x)dx\\
&+\e\int_{\R^n}\frac{\mu}{(1+|x|)^{\beta}}f(x)b_q(0, x)dx\Big)+C\int_0^T\int_{\R^n}u\Big[\partial_t^2(b_q\eta_M^{2p'})\\
&-\Delta(b_q\eta_M^{2p'})-\frac{\mu}{(1+|x|)^\beta}
\partial_t(b_q\eta_M^{2p'})\Big]dxdt\\
\le&C \int_0^T\int_{\R^n}u\left(2\partial_tb_q\partial_t\eta_M^{2p'}+b_q\partial_t^2\eta_M^{2p'}
-\frac{\mu}{(1+|x|)^\beta}b_q\partial_t\eta_M^{2p'}\right)\\
\triangleq& I_6+I_7+I_8.\\
\end{aligned}
\end{equation}

Next we will use the asymptotic behavior \eqref{bq3} to estimate $I_6-I_8$. First,
\begin{equation}\label{I6}
\begin{aligned}
I_6\le& CM^{-1}\left(\int_0^T\int_{\R^n}|u|^pb_q\theta_M^{2p'}dxdt\right)^{\frac1p}
\left(\int_0^T\int_{\R^n}b_q^{-\frac{1}{p-1}}b_{q+1}^{\frac{p}{p-1}}dxdt\right)^{\frac{p-1}p}\\
\le&CM^{-1}\left(\int_0^T\int_{\R^n}|u|^pb_q\theta_M^{2p'}dxdt\right)^{\frac1p}\\
&\times\left(\int_{\frac M2}^M\int_{0}^{t+1}(t+2+r)^{\frac{n-1}2-\frac1{p(p-1)}}(t+2-r)^{-1}drdt\right)^{\frac{p-1}p}\\
\le&C\left(\log M\right)^{\frac{p-1}p}\left(\int_0^T\int_{\R^n}|u|^pb_q\theta_M^{2p'}dxdt\right)^{\frac1p}.\\
\end{aligned}
\end{equation}
For $I_7$ and $I_8$, we have
\begin{equation}\label{I7}
\begin{aligned}
I_7\le& CM^{-2}\left(\int_0^T\int_{\R^n}|u|^pb_q\theta_M^{2p'}dxdt\right)^{\frac1p}
\left(\int_0^T\int_{\R^n}b_qdxdt\right)^{\frac{p-1}p}\\
\le&CM^{-2}\left(\int_0^T\int_{\R^n}|u|^pb_q\theta_M^{2p'}dxdt\right)^{\frac1p}\\
&\times\left(\int_{\frac M2}^M\int_{0}^{t+1}(t+2+r)^{-\left(\frac{n-1}2-\frac1{p}\right)}r^{n-1}drdt\right)^{\frac{p-1}p}\\
\le&C\left(\int_0^T\int_{\R^n}|u|^pb_q\theta_M^{2p'}dxdt\right)^{\frac1p}\\
\end{aligned}
\end{equation}
and
\begin{equation}\label{I8}
\begin{aligned}
I_8\le& CM^{-1}\left(\int_0^T\int_{\R^n}|u|^pb_q\theta_M^{2p'}dxdt\right)^{\frac1p}
\left(\int_0^T\int_{\R^n}\frac1{(1+|x|)^{\beta p'}}b_qdxdt\right)^{\frac{p-1}p}\\
\le&CM^{-1}\left(\int_0^T\int_{\R^n}|u|^pb_q\theta_M^{2p'}dxdt\right)^{\frac1p}\\
&\times\left(\int_{\frac M2}^M\int_{0}^{t+1}(t+2+r)^{-\left(\frac{n-1}2-\frac1{p}\right)}\frac{(1+r)^{n-1-p'}}
{(1+r)^{p'(\beta-1)}}drdt\right)^{\frac{p-1}p}\\
\le&C\left(\int_0^T\int_{\R^n}|u|^pb_q\theta_M^{2p'}dxdt\right)^{\frac1p}.\\
\end{aligned}
\end{equation}
Therefor we conclude from \eqref{Ycriineq}-\eqref{I8} that
\begin{equation}\label{Yp}
\begin{aligned}
Y^p(M)\le CM\left(\log M\right)^{p-1}Y'(M).
\end{aligned}
\end{equation}
The lifespan \eqref{lifespancri} yields by combining \eqref{Ycrideriv} and \eqref{Yp}, and using the following lemma with $p_1=p_2=p$
and $\delta=\e^p$, due to the fact $M$ is arbitrary in $(1, T)$.
\begin{lem}
{\bf (Lemma 3.10 in \cite{ISWa})}. Let $2<t_0<T$. $0\le \phi\in C^1([t_0, T))$. Assume that
\begin{equation}
\left\{
\begin{aligned}
& \delta\le K_1t\phi'(t), \quad t\in (t_0, T), \\
& \phi(t)^{p_1}\le K_2t(\log t)^{p_2-1}\phi'(t), \quad t\in (t_0, T)\\
\end{aligned}
\right.
\end{equation}
with $\delta, K_1, K_2>0$ and $p_1, p_2>1$. If $p_2<p_1+1$, then there exists positive constants $\delta_0$ and $K_3$(independent of $\delta$) such that
\begin{equation}
\begin{aligned}
T\le \exp\left(K_3\delta^{-\frac{p_1-1}{p_1-p_2+1}}\right)
\end{aligned}
\end{equation}
when $0<\delta<\delta_0$.
\end{lem}

\section{Proof for Theorem \ref{thmGlassey}}

Now we come to Theorem \ref{thmGlassey}. For $M\in (1, T)$, we set
\[
\psi(t, x)=-\eta_M^{2p'}(t)\Phi(t, x)=-\eta_M^{2p'}(t)e^{-t}\phi(x).\\
\]
It is interesting to see that
\begin{equation}\label{thmGlasseytest}
\begin{aligned}
\partial_t\psi(t, x)=\eta_M^{2p'}\Phi-2p'\eta_M^{2p'-1}\partial_t\eta_M\Phi\ge \eta_M^{2p'}\Phi>0,
\end{aligned}
\end{equation}
due to the fact that $\eta_M(t)$ is a non-increasing function. Also it is easy to check that
\begin{equation}\label{thmGlasseytestini}
\begin{aligned}
&\partial_t\psi(0, x)=\phi(x),\\
&\partial_t^2\psi(0, x)=-\phi(x).\\
\end{aligned}
\end{equation}

Noting \eqref{thmGlasseytestini}, if we set $\Psi(t, x)=\partial_t\psi(t, x)$ in the definition \eqref{weaksol} for \eqref{Glassey} and making integration by parts, then we come to
\begin{equation}\label{thm2step1}
\begin{aligned}
&\e \int_{\R^n}g(x)\phi(x)dx+\int_0^T\int_{\R^n}|u_t|^p\partial_t\psi dxdt\\
=&-\int_0^T\int_{\R^n}u_t\partial_t^2\psi dxdt+\int_0^T\int_{\R^n}\nabla u\cdot\nabla\psi_t dxdt\\
&+\int_0^T\int_{\R^n}\frac{\mu}{(1+|x|)^{\beta}}u_t\partial_t\psi dxdt£¬\\
\end{aligned}
\end{equation}
which implies
\begin{equation}\label{thm2step2}
\begin{aligned}
&\e \int_{\R^n}g(x)\phi(x)dx+\e\int_{\R^n}\left(1+\frac{\mu}{(1+|x|)^{\beta}}\right)f(x)\phi(x)dx\\
&+\int_0^T\int_{\R^n}|u_t|^p\partial_t\psi dxdt\\
=&\int_0^T\int_{\R^n}u_t\partial_t^2\left(\eta_M^{2p'}\Phi\right)dxdt-\int_0^T\int_{\R^n}u_t\Delta\left(\eta_M^{2p'}\Phi\right)dxdt\\
&-\int_0^T\int_{\R^n}\frac{\mu}{(1+|x|)^{\beta}}u_t\partial_t\left(\eta_M^{2p'}\Phi\right)dxdt\\
=&\int_0^T\int_{\R^n}u_t\left(\partial_t^2\eta_M^{2p'}\Phi+2\partial_t\eta_M^{2p'}\partial_t\Phi-\frac{\mu}{(1+|x|)^{\beta}}
\partial_t\eta_M^{2p'}\Phi\right)dxdt\\
&+\int_0^T\int_{\R^n}u_t\eta_M^{2p'}\left(\partial_t^2\Phi-\Delta\Phi-\frac{\mu}{(1+|x|)^{\beta}}\partial_t\Phi\right)dxdt\\
=&\int_0^T\int_{\R^n}u_t\left(\partial_t^2\eta_M^{2p'}\Phi+2\partial_t\eta_M^{2p'}\partial_t\Phi-\frac{\mu}{(1+|x|)^{\beta}}
\partial_t\eta_M^{2p'}\Phi\right)dxdt\\
\triangleq& I_9+I_{10}+I_{11}.\\
\end{aligned}
\end{equation}
By H\"{o}lder inequality and \eqref{testpro}, we have
\begin{equation}\label{VI}
\begin{aligned}
|I_9|\le& \int_0^T\int_{\R^n}|u_t\partial_t^2\eta_M^{2p'}\Phi|dxdt\\
\le& CM^{-2}\left(\int_0^T\int_{\R^n}\theta_M^{2p'}\Phi|u_t|^p dxdt\right)^{\frac1p}\left(\int_{\frac M2}^M\int_{\R^n}\Phi dxdt\right)^{\frac{1}{p'}}\\
\le& CM^{-2}\left(\int_0^T\int_{\R^n}\eta_M^{2p'}\Phi|u_t|^p dxdt\right)^{\frac1p}\left(\int_{\frac M2}^M(t+1)^{\frac{n-1}{2}}dt\right)^{\frac{1}{p'}}\\
\le& CM^{-1-\left(\frac{1}{p-1}-\frac{n-1}{2}\right)\frac{1}{p'}}\left(\int_0^T\int_{\R^n}\theta_M^{2p'}\Phi|u_t|^p dxdt\right)^{\frac1p}.\\
\end{aligned}
\end{equation}
\begin{equation}\label{VII}
\begin{aligned}
|I_{10}|\le& C\int_0^T\int_{\R^n}|u_t\partial_t\eta_M^{2p'}\Phi|dxdt\\
\le& CM^{-1}\left(\int_0^T\int_{\R^n}\theta_M^{2p'}\Phi|u_t|^p dxdt\right)^{\frac1p}\left(\int_{\frac M2}^M\int_{\R^n}\Phi dxdt\right)^{\frac{1}{p'}}\\
\le& CM^{-\left(\frac{1}{p-1}-\frac{n-1}{2}\right)\frac{1}{p'}}\left(\int_0^T\int_{\R^n}\theta_M^{2p'}\Phi|u_t|^p dxdt\right)^{\frac1p}\\
\end{aligned}
\end{equation}
and
\begin{equation}\label{VIII}
\begin{aligned}
|I_{11}|\le C|I_{10}|\le CM^{-\left(\frac{1}{p-1}-\frac{n-1}{2}\right)\frac{1}{p'}}\left(\int_0^T\int_{\R^n}\theta_M^{2p'}\Phi|u_t|^p dxdt\right)^{\frac1p}.\\
\end{aligned}
\end{equation}
By combing \eqref{thmGlasseytest}, \eqref{thm2step2}, \eqref{VI}, \eqref{VII} and \eqref{VIII}, we obtain
\begin{equation}\label{thm2step3}
\begin{aligned}
&C_2(f, g)\e +\int_0^T\int_{\R^n}\eta_M^{2p'}\Phi|u_t|^p dxdt\\
\le&CM^{-\left(\frac{1}{p-1}-\frac{n-1}{2}\right)\frac{1}{p'}}\left(\int_0^T\int_{\R^n}\theta_M^{2p'}\Phi|u_t|^p dxdt\right)^{\frac1p}.\\
\end{aligned}
\end{equation}

We consider the quantity
\[
Y\left[\Phi|\partial_tu|^p\right](M)=\int_1^M\left(\int_0^T\int_{\R^n}w(t, x)\theta_\sigma^{2p'}(t)dxdt\right)\sigma^{-1}d\sigma,
\]
then in the same way as in \eqref{Ycri} we get
\begin{equation}\label{Y}
\begin{aligned}
Y\left[\Phi|\partial_tu|^p\right](M)
\le C\log2\int_0^T\int_{\R^n}\eta_M^{2p'}\Phi|u_t|^pdxdt\\
\end{aligned}
\end{equation}
and
\begin{equation}\label{Yderiv}
\begin{aligned}
\frac{d}{dM}Y\left[\Phi|\partial_tu|^p\right](M)=M^{-1}\int_0^T\int_{\R^n}\theta_M^{2p'}\Phi|u_t|^pdxdt.\\
\end{aligned}
\end{equation}
For simplicity, we denote $Y(M)$ for $Y\left[\Phi|\partial_tu|^p\right](M)$. Then by combining \eqref{thm2step3}, \eqref{Y} and \eqref{Yderiv}, we know there exist positive constants $C_3, C_4$ such that
\begin{equation}\label{IneqforY}
\begin{aligned}
M^{-\left(\frac{1}{p-1}-\frac{n-1}{2}\right)(p-1)+1}Y'(M)\ge \left(C_3\e+C_4Y(M)\right)^p,
\end{aligned}
\end{equation}
which leads to
\begin{equation}
\label{lifespanM}
M\le
 \left\{
 \begin{array}{ll}
 C\e^{-\left(\frac{1}{p-1}-\frac{n-1}{2}\right)^{-1}} & \mbox{for}\ 1<p<p_G(n),\\
 \exp\left(C\e^{-(p-1)}\right) & \mbox{for}\ p=p_G(n).
 \end{array}
 \right.
\end{equation}
Since $M$ is arbitrary in $(1, T)$, we then obtain the lifespan estimate \eqref{lifespan2}.
\section*{Appendix}

We are left with the proof of Lemma \ref{testcri}. The proof is parallel to that of Lemma 3.1 in \cite{YZ06} and we write out the sketch. Given $\eta\in[0,1]$, suppose $H_0$ and $H$ be the fundamental solution of
$$\Delta u -\eta^2 u- u_t=0,\ \ \ \ \ \Delta u -\eta^2 u-\eta Vu- u_t=0$$
in $\mathbb{R}^n\times (0,\infty)$, respectively. If we set $H_0=e^{-t\eta^2}G_0$ and $H=e^{-t\eta^2}G$, then $G_0$ and $G$ are the fundamental solution of
$$\Delta u - u_t=0,\ \ \ \ \ \Delta u -\eta Vu- u_t=0$$
respectively. By Theorem 1.1(a) and Remark 1.1 in \cite{Zhang03}, there exists a positive constant
\begin{equation}\label{C_eta}
C(\eta):=e^{-c\sup_{x,y\in \mathbb{R}^n,t>0}N(\eta V,t)(x,y)}
\end{equation}
such that
\begin{equation}\label{a_G}
C(\eta)G_0(x,t;y,0)\leq G(x,t;y,0)\leq G_0(x,t;y,0)=\frac{C_n}{t^{n/2}}e^{-\frac{|x-y|^2}{4t}}
\end{equation}
for all $x,y\in\mathbb{R}^n$ and $t > 0$.

In the following we show that $C(\eta)$ has a uniform positive lower bound for $\eta\in[0,1]$.
Noting that $c$ in $C(\eta)$ depends only on spatial dimension $n$ and $N(V,t)(x,y)$ is defined as
\begin{eqnarray*}
N(V,t)(x,y)&=&\int_0^{t/2}\int\frac{e^{-|z-y+(\tau/t)(y-x)|^2/4\tau}}{\tau^{n/2}}|V(z)|dzd\tau\\
&&+\int_{t/2}^t\int\frac{e^{-|z-y+(\tau/t)(y-x)|^2/4(t-\tau)}}{(t-\tau)^{n/2}}|V(z)|dzd\tau.
\end{eqnarray*}
It has been proved in Proposition 2.1 in \cite{Zhang03} that
$$N(V,t)\leq C(\|V\|_p+\|V\|_q)$$
with $p>\frac n2$ and $q<\frac n2$. Hence, for $V(x)=\frac{\mu}{(1+|x|)^{\beta}}, \beta>2$, one can always choose
$p>\frac n2$ and $\frac n{2+\delta}<q<\frac n2$ such that
$$N(V,t)\leq C\left(\frac1{p\delta}+\frac1{(2+\delta)q-n}\right)<\infty.$$
Observing that $N(\eta V,t)=\eta N(V,t)$, it is easy to see that
$$c_0:=\inf_{0\leq\eta\leq1}\{C(\eta)\}=C(1)>e^{-cC\left(\frac1{p\delta}+\frac1{(2+\delta)q-n}\right)}>0.$$
Consequently, $C(\eta)$ in \eqref{a_G} can be replaced by $c_0$. Multiplying with $e^{-t\eta^2}$ on both sides, we find that the following global bounds hold uniformly for any $\eta\in[0,1]$
\begin{equation}\label{a1}
c_0H_0(x,t;y,0)\leq H(x,t;y,0)\leq H_0(x,t;y,0).
\end{equation}

The remaining steps are almost the same as that in \cite{YZ06}, we omit the details.

\section*{Acknowledgment}
\par\quad
This work is completed when the first author visited Illinois Institute of Technology, he would like to express his sincere thank to Prof. Chun Liu for his discussion and warm hospitality.

The first author is supported by NSF of Zhejiang Province(LY18A010008), Postdoctoral Research Foundation of China(2017M620128, 2018T110332), the Scientific Research Foundation of the First-Class Discipline of Zhejiang Province
(B)(201601), the CSC(201708330548). The second author is partially supported by NSF of Zhejiang Province(LY18A010023).


\bibliographystyle{plain}

\end{document}